\documentclass[letterpaper, 10 pt, conference]{ieeeconf}  

\IEEEoverridecommandlockouts                              
\usepackage[english]{babel}

\let\labelindent\relax 

\usepackage[colorinlistoftodos,bordercolor=orange,backgroundcolor=orange!20,linecolor=orange,textsize=scriptsize]{todonotes}

\usepackage{amsfonts,amssymb,mathrsfs,amsmath,amscd,epsfig,setspace,amstext,color,graphics,verbatim}
\usepackage{hyperref}
\usepackage{fixmath}
\usepackage{mathtools}
\usepackage{mathrsfs}

\usepackage{phdalgo}
\usepackage{comment}
\usepackage[english]{babel}
\usepackage[T1]{fontenc}
\usepackage{xargs}
\usepackage{amsmath}
\usepackage{amsfonts}
\usepackage{amssymb}

 \usepackage{amsthm}
  
\usepackage{float}
\usepackage{adjustbox}
\usepackage{bm}
\usepackage{bbm}
\usepackage{dsfont}
\makeatletter
\let\NAT@parse\undefined
\makeatother
\usepackage[numbers,sort&compress]{natbib}

\newcommand{\psg}{{\sc SG}}
\newcommand{\mpsg}{{\sc MPSG}}
\newcommand{\polmin}{\mathscr{S}}
\newcommand{\polmax}{\mathscr{T}}
\usepackage{enumitem}

\usepackage[framemethod=tikz]{mdframed}

\usepackage{color}
\definecolor{myblue}{rgb}{0.122, 0.435, 0.698}
\newmdenv[innerlinewidth=0.5pt, roundcorner=4pt,linecolor=myblue,innerleftmargin=6pt,
innerrightmargin=6pt,innertopmargin=6pt,innerbottommargin=6pt]{bluebox}
\definecolor{myred}{rgb}{0.8, 0.1, 0.1}
\newmdenv[innerlinewidth=0.5pt, roundcorner=4pt,linecolor=myred,innerleftmargin=6pt,
innerrightmargin=6pt,innertopmargin=6pt,innerbottommargin=6pt]{redbox}

\usepackage{stmaryrd}

\SetSymbolFont{stmry}{bold}{U}{stmry}{m}{n}

\usepackage{cleveref}

\newtheorem{theorem}{Theorem}
\newtheorem{corollary}{Corollary}
\newtheorem{proposition}{Proposition}
\newtheorem{lemma}{Lemma}
\newtheorem{assumption}{Assumption}
\theoremstyle{example}
\theoremstyle{definition}
\newtheorem{example}{Example}
\newtheorem{remark}{Remark}
\newtheorem{definition}{Definition}

\Crefname{corollary}{Cor.}{Cors.}
\Crefname{equation}{Eq.}{Eqs.}
\Crefname{figure}{Fig.}{Figs.}
\Crefname{tabular}{Tab.}{Tabs.}
\Crefname{table}{Tab.}{Tabs.}
\Crefname{theorem}{Thm.}{Thms.}
\Crefname{definition}{Def.}{Defs.}
\Crefname{section}{Sec.}{Secs.}
\Crefname{proposition}{Prop.}{Props.}
\Crefname{assumption}{Asm.}{Asms.}
\Crefname{example}{Ex.}{Exs.}
\Crefname{algocf}{Algorithm}{Algorithms}

\Crefname{appsec}{Appendix}{Appendices}

\Crefname{claim}{Claim}{Claims}

\usepackage{algorithm}
\usepackage{algpseudocode}

\renewcommand{\geq}{\geqslant}
\renewcommand{\leq}{\leqslant}

\graphicspath{{../images/},{../../../SIMU/plots/keep/CDC_Disag},{../../SIMU/plots/keep/smartGridComm/}}
\DeclareGraphicsExtensions{.pdf,.jpeg,.png,.svg}

\newcommand{\argmin}[1]{\underset{#1}{\mathrm{argmin }} }
\newcommand{\argmax}[1]{\underset{#1}{\mathrm{argmax }} }

\newcommand{\R}{\mathbb{R}}
\newcommand{\E}{\mathbb{E}}

\newcommand{\Id}{\mathrm{Id}}

\newcommand{\mc}{\mathcal}

\let\ifXXX\iffalse

\newcommand{\cw}{\operatorname{cw}}
\newcommand{\myphi}{\varphi}
\title{\LARGE \bf   Solving Ergodic Markov Decision Processes and Perfect Information Zero-sum Stochastic Games by Variance Reduced Deflated Value Iteration}

\author{
Marianne Akian, St{\'e}phane Gaubert, Zheng Qu, Omar Saadi\\
\textit{Accepted for publication in the proceedings of the CDC'2019 conference.}
\thanks{The authors were partially supported by PGMO, a joint program of EDF, Orange, Thales and FMJH (Fondation Math\'{e}matique Jacques Hadamard) and by IRS ICODE.}%
\thanks{M. Akian, S. Gaubert and O. Saadi are with INRIA and CMAP, \'{E}cole polytechnique. Address: CMAP, Ecole polytechnique, Route de Saclay, 91128 Palaiseau Cedex, France.   Email: {\tt\scriptsize marianne.akian@inria.fr}, {\tt\scriptsize stephane.gaubert@inria.fr}, {\tt\scriptsize omar.saadi@polytechnique.edu}.}%
\thanks{Z. Qu is with Dep. Math., The University of Hong Kong, Address: The University of Hong Kong, Room 419, Run Run Shaw Building Pokfulam Road, Hong Kong. Email: {\tt\scriptsize  zhengqu@hku.hk}.}%
}

\begin{document}

\maketitle

\begin{abstract}
  Recently, Sidford, Wang, Wu and Ye (2018) developed an algorithm combining variance reduction techniques with value iteration to solve discounted Markov decision processes. This algorithm has a sublinear complexity when the discount factor is fixed. Here, we extend this approach to mean-payoff problems, including both Markov decision processes and perfect information zero-sum stochastic games.
  We obtain sublinear complexity bounds, assuming there is a distinguished state which is accessible from all initial states and for all policies.
  Our method is based on a reduction from the mean payoff problem to the discounted problem by a Doob h-transform, combined with a deflation technique. The complexity analysis of this algorithm uses at the same time the techniques developed by Sidford et al. in the discounted case and non-linear spectral theory techniques (Collatz-Wielandt characterization of the eigenvalue).
\end{abstract}

\section{Introduction}
{\bf Context}.\/ 
Markov decision processes, and more generally zero-sum two player stochastic games, are classical models to study sequential 
problems under uncertainty \cite{puterman2014markov,NS03}. They appear
in various applications ranging from engineering sciences, finance, economy,
to health care or ecology.
The dynamic programming method allows one to reduce the infinite horizon
problem, in which players optimize a discounted payoff, to a fixed
point problem, involving an order preserving and contracting map,
called {\em Bellman} or {\em Shapley operator}. 
Value iteration and policy iteration~\cite{puterman2014markov}
are two fundamental dynamic programming methods. The execution
time of these two algorithms is generally super-linear in the size of the input,
and it is desirable to develop accelerated algorithms, for well structured
huge scale instances.

Algorithms based on Monte-Carlo simulations can lead to improved scalability.
In a recent progress, Sidford et al.~\cite{sidford2018variance} combined value iteration algorithm with sampling and variance reduction techniques. They obtained an algorithm for discounted infinite-horizon MDPs that, remarkably,
is sublinear in a certain relevant regime of the parameters. 

In the mean payoff problem, the discount factor tends to $1$ and the bounds of execution time for value iteration blow up which excludes
to pass to the limit in the algorithm of~\cite{sidford2018variance}.

{\bf Contribution.} In the present paper, we study
a class of two-player {\em mean payoff problems},
and we still obtain a sublinear complexity result.
This applies in particular to the mean payoff problems for MDPs,
which correspond to the special one-player case.

To do so, we exploit a general method, first introduced in a previous
work of two of the authors~\cite{akian2013policy}, 
allowing one to reduce a class 
of mean payoff problems to discounted problems.
In this way, complexity results concerning the discounted
problem can be transferred (with some additional work)
to the mean payoff case.

This reduction combines a scaling argument (a combinatorial version
of Doob's h-transform arising in the boundary theory
of Markov processes~\cite{dynkin}) and a deflation technique:
to a mean payoff
problem, we associate a discounted problem, with a {\em state-dependent}
discount rate (\Cref{ErgToDisc}).
Another key idea is the use of
{\em weighted sup-norms}, in order to obtain contraction
rates for this new Shapley operator. Then, we
solve the mean payoff problem by calling twice
a variant of the algorithm of Sidford's et al.~\cite{sidford2018variance}:
we call first this variant to compute the diagonal scaling involved in
our reduction, and then to solve the discounted game obtained
after the reduction. 
We also note that the present variant includes
an extension of the algorithm of~\cite{sidford2018variance}
for one player to the two-player
case, this is an easier matter---the main novelty here is rather
the reduction from the mean payoff problem to value iteration in the discounted case.

The subclass of problems to which our method applies requires
the existence of a distinguished state $c$ to which all other states have access, for all policies of the two players. The maximum of the hitting time
of this state, for all policies, appears in
our complexity bound. There are instances in which this maximal
hitting time is $O(1)$, and then we end up with a sublinear regime.

Weighted sup-norms were already used by Bertsekas
and Tsitsiklis
to obtain contraction
results for value iteration
in the case
of stochastic shortest path problems~\cite{bertsekas1991analysis}.
Our main results include \Cref{Clarke} and \Cref{coro-cw},
which characterize the best contraction rate, with respect
to all possible weighted sup-norms, as the Collatz-Wielandt number of a certain convex monotone positively homogeneous map which we call the ``Clarke recession function''.
This also implies that a contraction estimate previous computed in~\cite{akian2013policy}
is indeed optimal, if all the actions are ``useful'' in a natural sense.

{\bf Comparison with other approaches}. Gupta, Jain and Glynn have
recently developed a Monte-Carlo version of {\em relative value iteration} to solve mean payoff problems~\cite{jain}. The convergence analysis requires
the Bellman operator be a strict
contraction
in the ``span seminorm''. This is a demanding condition.
For instance, in the $0$-player case,
this requires the transition matrix to be primitive,
where our reduction holds in more general
circumstances (the uniqueness of the final class suffices).
Wang developed in \cite{wang2017primal} an algorithm for mean-payoff MDPs (one player), that has a sublinear bound. This algorithm depends on a
mixing time and on a parameter $\tau$ measuring the distance
between the invariant measures attached to different policies
(the mixing times of~\cite{wang2017primal} should not be confused
with the hitting times used here, the finiteness of the former
implies the finiteness of the latter, but not vice versa).
There are instances in which 
the distance $\tau$
is exponential in the input size, whereas the hitting time
is linear.

The paper is organized as follows. In Section~\ref{s-DynProg},
we recall basic notions about zero-sum games. In Section~\ref{sec-prelim},
we present the main techniques allowing the reduction from
the mean payoff case to the discounted case. In Section~\ref{sec-SSG},
we present an adaptation of the variance reduction algorithm of~\cite{sidford2018variance}, allowing us to handle the operators obtained after the
deflation and h-transform reduction. In Section~\ref{sec-ergodic},
we derive the sublinear bounds for classes of mean payoff problems.
Examples are presented in Section~\ref{sec-example}.
Most proofs are omitted owing to the space constraint.

\section{Dynamic programming operators}
\label{s-DynProg}
\subsection{Shapley operators of perfect information zero-sum stochastic games with general discount factor}
\label{sec-def-pol}
We refer the reader to~\cite{NS03}
for background on stochastic games. We next briefly
recall the main notions and properties.

A perfect information two-player zero-sum stochastic game with general discount (\psg)
is described by the following data. We consider a finite state space
$S:=\{1,\dots,n\}$. For all $i\in S$, $A_i$ is a finite set representing the possible actions of player MIN in state $i$, and $B_{i,a}$ is a finite set representing the possible actions of player MAX in state $i$, when player MIN just played action $a$.
We denote by $E:=\{(i,a,b)~|~ i\in S, a \in  A_i, b\in B_{i,a}\}$ the set of all admissible
triples state-actions. For all $(i,a,b)\in E$, $P^{ab}_i$ is an element of $\Delta(S)$ the set of probability measures on $S$; we shall identify $P^{ab}_i$ to a row vector in $\R^n$, writing $P^{ab}_i=(P^{ab}_{ij})_{j\in S}$ where $P^{ab}_{ij}$ is the transition probability to the next state $j$, given the current state $i$ and the actions taken $a\in A_i, b\in B_{i,a}$. For all $(i,a,b)\in E$, $r^{ab}_i$ is a reward (real number) that MIN pays to MAX, and $\gamma^{ab}_i$ (real nonnegative number) is a discount factor.
We define
\[ R:=\max_{(i,a,b)\in E} |r_i^{ab}| \in \R_{+},
\quad \Gamma:=\max_{(i,a,b)\in E} \gamma_i^{ab} \in (0,\infty)
\enspace,
\]
where $\R_{+}:= \{x\in \R\mid x\geq 0\}$.
We allow $\gamma_i^{ab}$ to take values larger than $1$.
The term
{\em turn-based} is sometimes used as a synonym of ``perfect information''. This is in contrast
with the more general model of Shapley's imperfect information stochastic games in which two players play simultaneously with
randomized actions, see e.g.~\cite{MN81}.

Recall that a {\em strategy} of a player is a decision rule which
associates to a history of the game an admissible action of this player.
A strategy $\sigma$
of player Min, a strategy $\tau$ of player Max,
and an initial state $i$, alltogether determine a
random process $(i_\ell, a_\ell,b_\ell)_{\ell\geq 0}$ with values
in $E$: $i_\ell$ represents the state at step $\ell$,
and $a_\ell,b_\ell$ represent the actions of the two
players at the same state. We require that $i_0=i$.
We denote by $\E_{i,\sigma,\tau}$ the expectation
with respect to the probability measure governing
this process. Given a finite horizon $k$, we consider the zero-sum game
in which the payoff of player Max is given by
\begin{align}
J_i^k(\sigma,\tau) = \E_{i,\sigma,\tau} \left(
\sum_{\ell=0}^{k-1} \big(\prod_{m=0}^{\ell-1} \gamma_{i_m}^{a_mb_m} )r_{i_l}^{a_\ell b_\ell}
\right)\enspace.
\label{e-generalfunctional}
\end{align}
The {\em value} $v^k_i$ of the $k$-stage game starting from $i$ is
\ifXXX\todo{add references}\fi
defined as
\begin{align}
v^k_{i}:= \inf_\sigma\sup_\tau J_i^k(\sigma,\tau)
=
\sup_\tau\inf_\sigma  J_i^k(\sigma,\tau) \enspace,\label{e-def-value}
\end{align}
where the infima and suprema are taken over the set
of strategies of both players. By definition, the existence of the value
requires the infimum and the supremum to commute. A pair
of strategies $\sigma^*,\tau^*$ is said to be {\em optimal}
if $\sigma^*$ achieves the first infimum in~\eqref{e-def-value}
and if $\tau^*$ achieves the second supremum in~\eqref{e-def-value}.

For all $(i,a,b)\in E$, we set $M^{ab}_{ij}:=\gamma_i^{ab} P^{ab}_{ij}$
and $M^{ab}_i:= (M_{ij}^{ab})_{j\in S} \in \R^n$.
\begin{definition}\label{def-shapley}
  For a given \psg, the Shapley operator $T$ is the map $\R^n \rightarrow \R^n$ whose $i$th coordinate is given by
	\begin{align*}
	  T_i(v)=\min_{a\in A_i} \max_{b\in B_{i,a}}\big\{r_i^{ab}+
          \sum_{j\in S}{M_{ij}^{ab}} v_j\big\} , \quad i\in S, \; v\in \R^n  \enspace .
	\end{align*}
\end{definition}
It is known that the value vector $v^k=(v^k_i)_{i\in S}$
does exist and satisfies
\begin{align}
  v^k = T(v^{k-1}) \enspace, \qquad v^0 =0 \enspace .
  \label{e-newdp}
\end{align}
Relations of this form are often established when the discount factor is constant~\cite{NS03}, they remain valid for general functionals of the form~\eqref{e-generalfunctional} with state and action dependent factors, see Chapter 11 of~\cite{whittle86}. Note that the assumption that the discount factor be smaller than $1$ is not needed for the well posedness of the finite horizon
problem and for the validity of~\eqref{e-newdp}.

Similarly, one can consider the {\em infinite
horizon discounted zero-sum game}, in which the payoff of player Max is now
\[
J_i(\sigma,\tau) = \E_{i,\sigma,\tau} \left(
\sum_{\ell=0}^{\infty} \big(\prod_{m=0}^{\ell-1} \gamma_{i_m}^{a_mb_m} )r_{i_\ell}^{a_\ell b_\ell}
\right) \enspace .
\]
This payment is well defined, in particular, when, $\Gamma<1$
since then the above series become absolutely convergent.
Then, the value of the infinite horizon game and the notion
of optimal strategies are defined
in a similar manner to the finite horizon case. The value
vector $v=(v_i)_{i\in S}$ does exist and it is characterized
as the unique solution of the fixed point problem
\[
v= T(v)  \enspace ,
\]
see again~\cite{NS03} for background.
We shall see later on that the assumption $\Gamma<1$ can be relaxed:
what matters is that the discount factor be smaller than one in an ``average'' sense.

In what follows, it will be convenient to consider a special
class of strategies, determined by {\em policies} (feedback, stationary
rules). A {\em policy} of player MIN is a map:
\begin{align*}
\sigma: S \rightarrow \cup_{i\in S}A_i \enspace ,\quad
i \mapsto \sigma(i)\in A_i \enspace .
\end{align*}
We denote by $\polmin$ the set of all policies of player MIN.
Similarly, a {\em policy} of MAX is a map:
\begin{align*}
  \tau: \cup_{i\in S}(i,A_i) &\rightarrow \cup_{i\in S,\, a\in A_i }B_{i,a} \enspace ,
(i,a) \mapsto \tau(i,a)\in B_{i,a} \enspace .
\end{align*}
Note that since the game is in perfect information, MAX observes the action $a$ of MIN, and so the policy of player MAX takes care of this action.
We denote by $\polmax$ the set of all policies of player MAX.
It is known that in the discounted game, there exist optimal strategies
associated to policies (the action is selected at each step
by applying a policy of one player, the policy being the same
for all time steps). These policies are obtained by selecting actions
achieving the minimum and the maximum in the expression of $T(v)$
in \Cref{def-shapley}. See~\cite{NS03}.

Any choice of policies $(\sigma,\tau)\in \polmin\times \polmax$ defines the Markovian matrix $P^{\sigma\tau}\in \R^{n\times n}$ which determines the state trajectory
if the two players select their actions according to these policies.
i.e., $(P^{\sigma\tau})_{ij}=P^{\sigma(i) \tau(i,\sigma(i))}_{ij}$.
Similarly, we define the nonnegative matrix
$M^{\sigma\tau}$ with entries  $(M^{\sigma\tau})_{ij}=M^{\sigma(i) \tau(i,\sigma(i))}_{ij}$.
We denote the cardinality of a finite set $\mc{S}$ by $|\mc{S}|$. We recall that the size of the input is of order $|S\|E|$.

\subsection{Mean payoff problem}
\label{MPP}
We are now interested in the {\em undiscounted} case, in which the discount
factor $\gamma$ is identically $1$. Then we are considering a two-player perfect information zero-sum Stochastic Mean-Payoff Game (\mpsg), where the main quantity of interest
is the {\em mean payoff vector}:
\[
\chi(T) := \lim_{k\to \infty} T^k(0)/k  \enspace .
\]
The entry $\chi_i(T)$ represents the mean payoff per time unit,
if the initial state is $i$. Here, the mean payoff is
defined by considering a family of games in finite horizon $k$
as $k$ tends to infinity. There are alternative approaches,
in which the mean payoff is defined as the value of an infinite horizon game~\cite{LL69}. The property of the {\em uniform value} established in~\cite{MN81}
entails that the different natural approaches lead to the same notion of mean payoff.

The analysis of the mean payoff
problem is simplified when the following {\em non-linear eigenproblem}
has a solution:
\begin{align}
  \eta e +v=T(v), \qquad \eta\in \R, \;v\in \R^n \enspace ,
  \label{e-ergodic}
\end{align}
where $e:=(1 \cdots 1)^\top \in \R^n$ is the unit vector.
The scalar $\eta$ is called the {\em ergodic constant}, whereas
the vector $v$, which is not unique, is called {\em bias}
or {\em  potential}.
When this equation is solvable, we have $\chi(T) = \eta e$, i.e., the mean payoff is
independent of the initial state, and it is equal to the ergodic constant.
See e.g.~\cite{1610.09651} for background.

\section{Reduction from a mean payoff problem to a discounted one}\label{sec-prelim}

Let $u\in\R^n$, we write $u\gg0$ and we say that $u$ is a {\em positive vector} if for all $i\in [n]:=\{1,\cdots,n\}$, $u_i>0$.
Given $u\gg0$, we define the {\em weighted sup norm} $\|\cdot\|_u$ by :
\begin{equation*}
\|x\|_u=\max_{1\leq i\leq n} \frac{x_i}{u_i} =\|u^{-1}x\|_\infty \text{ , } \forall x\in \R^n \enspace , 
\end{equation*}
where the notation $u^{-1}x:=(u_i^{-1}x_i)_{i\in [n]}$ refers to the Hadamard quotient.
For $x$, $y$ $\in\R^n$ we write $x\leq y$ if $x_i\leq y_i$ for all $i\in[n]$.
A function $f:\R^n\rightarrow \R^n$ is said to be {\em monotone} if for all $x$, $y$ $\in\R^n$, if $x\leq y$ then $f(x) \leq f(y)$.
\subsection{Contraction rate of Shapley operators}\label{subsec-contraction}
We next introduce a notion of recession function associated
to a non-linear map.
Our definition is inspired by the 
notion of Clarke generalized directional derivative~\cite[Ch.~2, S1]{clarke} of a function $f$
at point $z$ in the direction $y$
\begin{align}
f'_z(y) := \limsup_{x\to z, \, s \to 0^+} \frac{f(x+s y) - f(x)}{s}
\enspace.\label{e-clarkeder}
\end{align}
We next adapt this idea by considering ``variations at infinity''
instead of local variations.
\begin{definition}
	Given a function $f:\R^n\rightarrow \R^n$, we define $\hat{f}:\R^n\rightarrow (\R\cup \{+\infty\})^n$ the {\em Clarke recession function} of $f$ as:
\begin{align}
  \hat{f}(y)=
  \sup_{s>0,\,  x\in \R^n}
  \frac{f(x+s y)-f(x)}{s} \enspace .
        \label{e-clarkerec}
\end{align}
\end{definition}
We chose the name ``Clarke recession function'' in view
of the similarity between~\eqref{e-clarkerec} and~\eqref{e-clarkeder}.

The following result is immediate:

\begin{proposition}
  The Clarke recession function is positively homogeneous and convex.
\end{proposition}
\begin{theorem}\label{Clarke}
  Let $f:\R^n\rightarrow \R^n$ be a monotone function, $u\gg0$ be a positive vector, and $\lambda \in \R_{+}$. We have $\hat{f}(u)\leq \lambda u$ if and only if the function $f$ is $\lambda$-contracting in the weighted sup-norm $\|\cdot\|_u$:
  \[\forall x,y \in\R^n,\qquad \|f(x)-f(y)\|_u\leq \lambda\|x-y\|_u
  \enspace .\]
\end{theorem}

Following~\cite{mallet2002eigenvalues,akian2013policy}, we define the
{\em Collatz-Wielandt number} of $\hat{f}$ as
\[\cw(\hat{f}):=\inf\{\lambda>0\mid  \exists u \gg0; \hat{f}(u)\leq \lambda u\}
\enspace .
\]
As an immediate consequence of \Cref{Clarke}, we get:
\begin{corollary}\label{coro-cw}
	If $f:\R^n\rightarrow \R^n$ is monotone, then:
	\begin{equation*}
	\cw(\hat{f})=\inf\{\lambda\in \R_{+} ~|~ \exists u \gg0; f \text{ is } \lambda\text{-contracting in } \|\cdot\|_u \} .
	\end{equation*}
\end{corollary}
We consider the Shapley operator $T:\R^n \rightarrow \R^n$
of \Cref{def-shapley}.
The following "max-max" operator $T^{\max}:\R^n \rightarrow \R^n$
was considered in
\cite{akian2013policy}
\begin{equation}
T^{\max}_i(y)=\max_{a\in A_i,b\in B_{i,a}} \{M_i^{ab} y\} \text{ , } \forall i\in S, y\in \R^n \enspace.
\end{equation}
In contrast to the Clarke recession function $\hat{T}$, 
$T^{\max}$ generally depends on the choice of the representation
of $T$ as a minimax expression. We next show, however,
that $T^{\max} = \hat{T}$ if all the terms arising
in the minimax expression are ``useful'' in the following
sense.
\begin{definition}
  For a given couple of actions $(a,b)$ of the two players, we define the set $C^{ab}_i=\{x\in \R^n \mid T_i(x)=r_i^{ab}+M_i^{ab} x\}$. We say that the couple of actions $(a,b)$ is {\em useful}
    if $\operatorname{int}(C^{ab}_i)\neq \emptyset$ for all $i\in [n]$.
\end{definition}
In the one player case, checking whether one action is useful 
reduces to checking whether a polyhedron has a non-empty
interior, and this can be done in polynomial time.

\begin{lemma}\label{Tmax}
	The Clarke recession function of the Shapley operator $T$ satisfies the following inequality:
	\begin{equation}
	\hat{T}(y)\leq T^{\max}(y) \text{ , } \forall y\in \R^n \enspace .
	\end{equation}
        Moreover, the equality holds if all the actions $(a,b)$ are useful.
\end{lemma}

There is an explicit formula for the Collatz-Wielandt number of $T^{\max}$.
  Recall that the notation $M^{\sigma\tau}$ refers to the nonnegative matrix
  associated to a pair of policies (end of Section~\ref{sec-def-pol}).
  We denote by $\rho(\cdot)$ the spectral radius of a matrix.
  \ifXXX\todo[inline]{SG: more precise ref}\fi
\begin{theorem}[{Corollary of~\cite{akian2011collatz}}]
  We have
  \[
  \cw(T^{\max}) = \max_{\sigma\in \polmin,\tau\in \polmax}\rho (M^{\sigma\tau})
  \enspace .\]
\end{theorem}

Owing to \Cref{Clarke} and \Cref{Tmax}, we will look for a vector $\myphi\gg0$ such that $T^{\max}(\myphi)\leq \lambda \myphi$ for some $\lambda \in [0,1)$, to have that the Shapley operator $T$ is $\lambda-$contracting in the weighted norm $\|\cdot\|_\myphi$. The following special construction allows us to obtain
such a $\myphi$ by solving a non-linear eigenproblem.

\begin{theorem}[Th.~7 and proof of Th.~13 in {\cite{akian2013policy}}]\label{PhiContraction}
  The following assertions are equivalent:
  \begin{enumerate}
    \item $\max_{\sigma\in \polmin,\tau\in \polmax}\rho (M^{\sigma\tau})<1$;
  \item there exists a unique vector $\myphi\in\R_{+}^n$ such that $\myphi=e+T^{\max}(\myphi)$.
  \end{enumerate}
  When these assumptions are satisfied, $T$ is $\lambda$-contracting in the weighted norm $\|\cdot\|_\myphi$, with $\lambda:=1-\frac{1}{\|\myphi\|_\infty}$.
\end{theorem}

\subsection{Correspondence between ergodic problems and discounted problems via the h-transform}\label{subsec-doob}
We consider here the 
non-linear eigenproblem~\eqref{e-ergodic},
where $T$ is the Shapley operator in the undiscounted case,
and describe a technique introduced in~\cite{akian2013policy}
to reduce this equation to a fixed point equation of a contracting operator.
Recall that~\eqref{e-ergodic} allows one to solve the  mean-payoff problem.
As noted above, the vector $v$ solution of~\eqref{e-ergodic} 
is not unique. In particular, if $v$ is a solution then $v + \alpha e$
also yields a solution for all $\alpha\in \R$. 
Hence, we shall distinguish
a special state $c\in S$ and require $v_c=0$.

\begin{definition}
	For a Markov matrix $P$ and states $i, j$, we denote:
	\begin{equation*}
	\mc{T}_{ij}(P):=\mathbb{E}[\text{inf}\{k\geq 1 ~|~ X_k=j\} ~|~ X_0=i]
	\end{equation*}
	the expected first hitting time of state $j$,
	for a Markov chain $X_k$ with transition matrix $P$ and initial state $i$.
\end{definition}
Given $c\in S$, it is easy to see that $\mathcal{T}_{ic}(P)<+\infty$ for all $i\in S$ if and only if $P$ has a unique final (recurrent) class and that $c$ belongs to this class. A state $c$ with the latter property is called a {\em renewal} state.

\begin{definition}
  For any matrix $P\in \R^{n\times n}$, we denote by $P_{(c)}\in \R^{n\times n}$ the matrix obtained from $P$ by replacing the column $c$ of $P$ with zeros.
  We denote by $P_i$ the $ith$ row of $P$, so that $P_i=(P_{ij})_{j\in[n]}$,
  and we use a similar notation for matrices constructed from $P$, e.g.,
  $P_{(c)i}=((P_{(c)})_{ij})_{j\in [n]}$. 
\end{definition}

\begin{lemma}\label{RhoHypothesis}
  Let $c\in S$ be a given state. The following assertions
  are equivalent:
  \begin{enumerate}
    \item	For all $(\sigma,\tau)\in \polmin\times \polmax$, $P^{\sigma\tau}$
    has a unique final class, and the state $c$
    is common to each of these classes;
	\item
	$\mathcal{T}_{ic}:=\max_{\sigma\in \polmin,\tau\in \polmax}\mc{T}_{ic}(P^{\sigma\tau})<+\infty, \quad \forall i \in S$;
	\item $\max_{\sigma\in \polmin,\tau\in \polmax}\rho (P^{\sigma\tau}_{(c)})<1$.
	\item
   There is a unique vector $\myphi^\star$ solution of the equation:
	\begin{equation}\label{PhiTmax}
	\myphi^\star=e+\max_{\sigma\in \polmin,\tau\in \polmax}[P^{\sigma\tau}_{(c)}\myphi^\star]\enspace,
	\end{equation}
        \end{enumerate}
        Under these assumptions, we have $\myphi^\star_i=\mc{T}_{ic}$, for all $i \in S$.
\end{lemma}
In the rest of this section, we make the following assumption.

\begin{assumption}\label{assum-commonc}
  There exists a state $c\in S$ satisfying the conditions
  of \Cref{RhoHypothesis}.
\end{assumption}
We can find such a state $c$ if it exists, or certify that there is none, in quasi-linear time by using directed hypergraphs techniques, along the lines of~\cite{1405.4658};
details will be given elsewhere.

Let $\myphi \in \R_{+}^n$, $\myphi\gg0$, and
$\R^n_c:= \{x\in \R^n\mid x_c=0\}$.
\begin{definition}
	For a nonnegative matrix $P \in \R^{n\times n}$, if
	$\myphi_i\geq 1+P_{(c)i}\myphi, \forall i\in S$, then we denote by $P_{(c,\myphi)}$ the nonnegative matrix obtained from $P$ by replacing the column $c$ by the vector $\myphi^{-1}_c(\myphi-1-P_{(c)}\myphi)$.
\end{definition}

\begin{lemma}\label{cPhi}
	Let $\eta\in \R$, $v\in \R^{n}$ with $v_c=0$ and $P \in \R^{n\times n}$.
	We have $\eta(\myphi-1)+Pv=P_{(c,\myphi)}(\eta \myphi+v)$.
	In particular $P_{(c,\myphi)}\myphi=\myphi-1$.
\end{lemma}

\begin{proposition}
The map $L_\myphi: (\eta,v) \mapsto w=\eta + \myphi^{-1}v$
from $\R\times \R^n_c$ to $\R^n$,
is an isomorphism, with inverse given by $w\mapsto (\eta,v)$ with
$\eta=w_c$ and $v=\myphi(w-w_c)$.
\end{proposition}
\begin{definition}
For any self-map $f$ of $\R^n$, we denote by
$\mathcal{L}_{\myphi}(f)$ the self-map of $\R^n$,
such that for all $w,v \in \R^n$ and $\eta \in \R$ with $v_c=0$ and $w=\eta + \myphi^{-1}v$, we have
\[
\mathcal{L}_{\myphi}(f)(w)=\myphi^{-1}(\eta(\myphi-1)+f(v))\enspace.
\]	
\end{definition}
\begin{remark}
	For a matrix $P$ we have by \Cref{cPhi} $\mathcal{L}_{\myphi}(P)(w)=P_{\myphi}w$ where $P_{\myphi}\in \R^{n\times n}$ is the nonnegative matrix given by $P_{\myphi,ij}:= \myphi^{-1}_iP_{(c,\myphi)ij} \myphi_{j}, (i,j)\in S^2$, so that
	\begin{align*}
	P_{\myphi,ij}=\left\{ \begin{array}{ll}
	\myphi^{-1}_i P_{ij} \myphi_j, \text{ if } j\neq c,\; i\in S \\
	1-\myphi^{-1}_i-\sum_{k\neq c}\myphi^{-1}_i P_{ik} \myphi_k, \text{ if } j=c,\; i\in S
	\end{array} \right..
	\end{align*}
\end{remark}
We consider the Shapley operator in the undiscounted case
\begin{equation}
T_i(v)=\min_{a\in A_i} \max_{b\in B_{i,a}}\big\{r_i^{ab}+P_{i}^{ab} v\big\} , \quad \forall i\in S, \forall v\in \R^n  .
\end{equation}
By \Cref{RhoHypothesis}, we know that there exists a vector $\myphi \in \R_{+}^n$, such that 
\begin{equation}\label{PhiIneq}
\myphi_i\geq 1+\max_{a,b}[P^{ab}_{(c)i}\myphi], \forall i\in S\enspace.
\end{equation}
So we can define as above a monotone operator $T^{\myphi}:=\mathcal{L}_{\myphi}(T)$, and we verify easily that
\begin{equation}\label{TphiExpression}
T_i^{\myphi}(w)= \min_{a\in A_i} \max_{b\in B_{i,a}}\big\{\myphi_i^{-1}r_i^{ab}+P_{\myphi, i}^{ab} w \big\}\enspace.
\end{equation}
\begin{lemma}\label{TPhiContraction}
	If there exists a vector $\myphi \in \R_{+}^n$ such that $\myphi_i\geq 1+\max_{a,b}[P^{ab}_{(c)i}\myphi]$ for all $i\in S$, then $T^{\myphi}$ is $\lambda_\myphi-$contracting in the sup-norm $\|\cdot\|_\infty$, with $\lambda_\myphi:=1-1/\|\myphi\|_\infty$. $T^{\myphi}$ can be interpreted as a Shapley operator of a discounted game with discount factors $\leq \lambda_\myphi$.
\end{lemma}

\begin{remark}\label{Tphimax}
  The vector $\myphi^\star$ defined by~\eqref{PhiTmax} is
  solution of a fixed point equation of the form $\myphi^\star=F(\myphi^\star)$
  where the map $F$ is order preserving and contracting in the norm $\|\cdot\|_{\myphi^\star}$. It follows that 
  if $w \in \R_{+}^n$, verifies $w_i\geq 1+\max_{a,b}[P^{ab}_{(c)i}w], \forall i\in S$, then, 
  $w\geq \myphi^\star$. Similarly,
  if $w_i\leq 1+\max_{a,b}P^{ab}_{(c,\myphi)i}w,\forall i\in S$ then $w\leq \myphi^\star$.
\end{remark}

\begin{theorem}\label{ErgToDisc}
  The non-linear eigenproblem 
	\begin{equation}
\label{eq-erg}
	\eta e+v=T(v)\enspace,
	\end{equation}
	where $\eta\in \R$ and $v\in \R^n$ with $v_c=0$, can be reduced to the fixed point problem:
	\begin{equation}\label{Tphi}
	T^{\myphi}(w)=w \enspace,
	\end{equation}
	where $w\in \R^n$ such that $w=\eta + \myphi^{-1}v$. Equation~\eqref{Tphi} has a unique solution $w^\star$.
\end{theorem}

We verify easily that for $w\in \R^n$, and $i\in S$ we have 
\begin{equation}\label{T-Phi}
\begin{split}
T_i^{\myphi}(w)= \min_{a\in A_i} \max_{b\in B_{i,a}}\big\{\myphi^{-1}_i P^{ab}_{i} \myphi(w-w_ce)\\
+\myphi^{-1}_i r^{ab}_{i}+w_c(1-\myphi_i^{-1}) \big\} , \quad \forall i\in S.
\end{split}
\end{equation}

\begin{lemma}\label{boundPay}
	The solution $w^\star$ of equation \eqref{Tphi} satisfies $\|w^\star\|_\infty\leq R$.
\end{lemma}

\begin{example}\label{ex-illustrate}
We give an elementary illustration
of the present
deflation+h-transform technique.
Let $P=\left(\begin{smallmatrix}0&1\\1&0\end{smallmatrix}\right)$, $r\in \R^2$, and consider $T(x)=r+Px$.
Let us choose $c=1$. The first hitting time vector $\myphi^\star$ is such that $\myphi^\star_2=1$ and $\myphi^\star_1=1+\myphi^\star_2=2$, so $\|\myphi^\star\|_\infty=2$. The operator $T^{\myphi^\star}$
given by~\eqref{T-Phi} specializes to
\(
T^{\varphi}(w_1,w_2)= ( \frac{r_1}{2} + \frac{w_2}{2},  r_2)  
\).
In accordance with~\Cref{TPhiContraction}, this operator is $1/2-$contracting.
The unique fixed point of $T^{\varphi}$ is $w=((r_1+r_2)/2, r_2)= \eta e + \myphi^{-1}v$, from which, by \Cref{ErgToDisc},
we recover the mean payoff $\eta=(r_1+r_2)/2$,
and $v=(0, (r_2-r_1)/2)$.
\end{example}
\section{Variance reduced value iteration for structured stochastic games}
\label{sec-SSG}

To solve the non-linear eigenproblem~\eqref{eq-erg},
we will find a vector $\myphi$ satisfying
\eqref{PhiIneq} by solving~\eqref{PhiTmax} in an approximate way,
and then use $\myphi$ to define
the new operator \eqref{TphiExpression} and solve 
the discounted problem \eqref{Tphi}.
We next present a variant of the method of Sidford et al.~\cite{sidford2018variance}
to deal with a structured
input, which will allow us to handle both problems
\eqref{PhiTmax} and \eqref{Tphi}.

We consider a perfect information two-player zero-sum stochastic game with general discount (\psg) as described in \cref{s-DynProg}, except that we suppose that $P^{\sigma \tau}$ is a sub-Markovian matrix for each couple of policies $(\sigma,\tau)\in \polmin\times\polmax$.
We suppose the associated Shapley operator can be written as
\begin{equation}\label{Tstructured}
T_i(w)= \min_{a\in A_i} \max_{b\in B_{i,a}}\big\{\gamma^{ab}_{i} P^{ab}_{i} Lw+G^{ab}_{i}(w) \big\}, \forall i\in S .  
\end{equation}
Here $L\in \R^{n\times n}$ is a sparse operator such that for all $w\in \R^n$, $Lw$ can be computed in $O(|S|)$. For all $i\in S, a\in A_i, b\in B_{i,a}$, $G^{ab}_{i}$ is a sparse affine operator such that $G^{ab}_{i}(w)$ can be computed in $O(1)$ for all $w\in \R^n$.
For example, by taking $L=\Id$ and $G^{ab}_{i}(w)=r^{ab}_{i}, \forall w\in\R^n$ we obtain the Shapley operator of the stochastic game with general discount.
The operator $L$ will allow us to handle the deflation (pre-subtraction of $w_c e$) in~\Cref{T-Phi}.

The problem that we want to solve is:
\begin{equation}\label{T}
T(w)=w \enspace.
\end{equation}

In this section, we make the following
assumption:
\begin{assumption}
\begin{itemize}
	\item $T$ is $\lambda-$contracting under the weighted norm $\|\cdot\|_\psi$, where $\psi\in \R^n$ is a positive vector.
	\item The solution $w^\star$ of the equation \eqref{T} verifies $\|w^\star\|_\psi\leq W$, where $W\geq 0$ is a scalar.
\end{itemize}
\end{assumption}

We can easily show the following inequalities: 
\begin{equation*}\label{NormIneq}
\|\psi^{-1}\|_\infty^{-1} \|w\|_\psi\leq \|w\|_\infty \leq \|\psi\|_\infty \|w\|_\psi \; \forall w\in \R^n\enspace.
\end{equation*}
\begin{remark}\label{d1d2}
	In the following, the values $\|\psi^{-1}\|_\infty$ and $\|\psi\|_\infty$ can be replaced by any positive scalars $d_1,d_2>0$ such that $\|\psi^{-1}\|_\infty\leq d_1$ and $\|\psi\|_\infty\leq d_2$.
\end{remark}

We denote also by $\|\cdot\|_\infty$ the operator norm associated to the sup-norm, so that we have:
\[
\|Mw\|_\infty\leq \|M\|_\infty\|w\|_\infty, \forall w\in \R^n, \forall M\in \R^{n\times n}.
\]

To a given vector $p=(p_j)_{j\in S}$ with $p_j\geq 0, \forall j\in S$ and $\sum_{j\in S}p_j\leq 1$, we associate the probability vector $\bar{p}=(\bar{p}_j)_{j\in S\cup \{0\}}$ with $\bar{p}_j=p_j, \forall j\in S$ and $\bar{p}_0=1-\sum_{j\in S}p_{j}$.

For each $i\in S$, $a\in A_i$, $b\in B_{i,a}$, we suppose that we can sample under the probability $\bar{P}_i^{ab}$ associated to the vector $P_i^{ab}$ in time $O(1)$.

We next adapt the algorithms $1-6$ presented by Sidford et al. in \cite{sidford2018variance} to our case with two players.
We follow the presentation of~\cite{sidford2018variance},
including the decomposition of the algorithm
in elementary subroutines. The necessary changes arise from
the use of the weighted sup-norm $\|\cdot\|_\psi$ instead of $\|\cdot\|_\infty$,
from the sub-Markovian character of the matrices. 

In the following, \Cref{Algo1} computes an approximation of $P^{ab}_i u$ by sampling under the probability vector $\bar{P}_i^{ab}$. \Cref{Algo2} computes an approximation of $T(w)$ for $w\in \R^n$, given an initial vector $w_0\in \R^n$. This algorithm assumes that
an approximation of the terms $x_{i}^{ab}=P_i^{ab}Lw_0$, called offsets in~\cite{sidford2018variance}, is already known.
Then, \Cref{Algo3} implements a randomized value iteration, using \Cref{Algo2} at each iteration.  To initialize \Cref{Algo3}, the
offsets $x_{i}^{ab}=P_i^{ab}Lw_0$ are computed exactly. \Cref{Algo4} iterates \Cref{Algo3}, using the technique of variance reduction by dividing the error by $2$ at every iteration. \Cref{Algo5} and \Cref{Algo6} are similar
to \Cref{Algo3} and \Cref{Algo4}, with the difference that the offsets are sampled, instead of being computed exactly.

\begin{algorithm}[h]
	\begin{small}
		\begin{algorithmic}[1]
			\\ \Lcomment{Input: vector $u\in \R^n$ and $M\geq 0$ such that we have $\|u\|_\infty\leq M$}
			\\ \Lcomment{Input: State $i\in S$ and actions $a\in A_i$, $b\in B_{i,a}$}
			\\ \Lcomment{Input: Target accuracy $\epsilon>0$, failure probability $\delta \in(0,1)$}
			\State $u_0=0$
			\State $m=\lceil\frac{2 M^2}{\epsilon^2}\ln(\frac{2}{\delta})\rceil$
			\For{$k\in [m]$} choose $i_k\in S\cup \{0\}$ with probabilities $\mathbb{P}(i_k=j)=\bar{P}^{ab}_{ij}$ for $j\in S\cup \{0\}$.
			\EndFor
			\\ \Return $Y=\frac{1}{m}\sum_{k\in [m]}u_{i_k}$
		\end{algorithmic}
	\end{small}
	\caption{Approximate transition with cemetery: \texttt{ApxTransC}{$(u,M,i,a,b,\epsilon,\delta)$}} \label{Algo1}
\end{algorithm}
\begin{algorithm}[h]
	\begin{small}
		\begin{algorithmic}[1]
			\\ \Lcomment{Input: Current vector $w \in \R^n$ and initial vector $w_0\in \R^n.$}
			\\ \Lcomment{Input: Precomputed offsets: $x \in \R^{E}$ with $|x^{ab}_i-P^{ab}_i Lw_0|\leq \epsilon$ for all $i\in S, a\in A_i$, $b\in B_{i,a}$.}
			\\ \Lcomment{Input: Target accuracy $\epsilon>0$, failure probability $\delta \in(0,1)$}
			\State  $M=\|L\|_\infty\|w-w_0\|_\infty$
			\State  $u=L(w-w_0)$
			\For{$i\in S$} 
			\For{$a\in A_i$}
			\For{$b\in B_{i,a}$}
			\State
			$\tilde{S}^{ab}_i=x^{ab}_i+\texttt{ApxTransC}(u,M,i,a,b,\epsilon,\frac{\delta}{|E|})$
			\State $\tilde{Q}^{ab}_i=\gamma^{ab}_i \tilde{S}^{ab}_i+G^{ab}_i(w)$
			\EndFor
			\State $\tilde{w}^a_i=\max_{b\in B_{i,a}} \tilde{Q}^{ab}_i$, $\tau(i,a) \in \argmax{b\in B_{i,a}}~\tilde{Q}^{ab}_i$
			\EndFor
			\State $\tilde{w}_i=\min_{a\in A_i}\tilde{w}^a_i$, $\sigma(i) \in \argmin{a\in A_i}~\tilde{w}^a_i$
			\EndFor
			\\ \Return $(\tilde{w},\sigma,\tau)$
		\end{algorithmic}
	\end{small}
	\caption{Structured approximate value operator: \texttt{SApxVal}{$(w,w_0,x,\epsilon,\delta)$}}
	\label{Algo2}
\end{algorithm}

\begin{algorithm}[h]
\begin{small}
	\begin{algorithmic}[1]
		\\ \Lcomment{Input: initial vector $w_0\in \R^n$, number of iterations $J>0$}
		\\ \Lcomment{Input: Target accuracy $\epsilon>0$, failure probability $\delta \in(0,1)$}
		\State \label{offsets} Compute $x\in \R^{E}$ such that $x^{ab}_i=P^{ab}_i Lw_0$ for all $i\in S$ and $a\in A_i$, $b\in B_{i,a}$.
		\For{$j\in [J]$} $(w_j,\sigma_j,\tau_j)=\texttt{SApxVal} (w_{j-1},w_0,x,\epsilon, \frac{\delta}{J})$
		\EndFor
		\\ \Return $(w_J,\sigma_J,\tau_J)$
	\end{algorithmic}
\end{small}
\caption{Structured randomized value iteration: \texttt{SRandVI}{$(w_0,J,\epsilon,\delta)$}} \label{Algo3}
\end{algorithm}

\begin{algorithm}[h]
	\begin{small}
		\begin{algorithmic}[1]
			\\ \Lcomment{Input: Target accuracy $\epsilon>0$, failure probability $\delta \in(0,1)$}
			\State Let $K=\lceil\log_2(\frac{\|\psi\|_\infty W}{\epsilon})\rceil$ and $J= \lceil\frac{1}{1-\lambda}\log(4)\rceil$
			\State $w_0=0$ and $\epsilon_0=W$
			\For{$k\in[K]$} \State $\epsilon_k=\frac {\epsilon_{k-1}}{2}=\frac{W}{2^k}$
			\State $(w_k,\sigma_k,\tau_k)=\texttt {SRandVI}(w_{k-1},J,\frac{1-\lambda}{4 \|\psi^{-1}\|_\infty \Gamma}\epsilon_k,\delta/K)$
			\EndFor
			\\ \Return $(w_K,\sigma_K,\tau_K)$
		\end{algorithmic}
	\end{small}
	\caption{Structured high precision randomized value iteration:\protect\\ \texttt{SHighPrecisionRandVI}{$(\epsilon,\delta,\lambda,W,\Gamma,\|\psi^{-1}\|_\infty,\|\psi\|_\infty)$}}
	\label{Algo4}
\end{algorithm}

\begin{algorithm}[h]
	\begin{small}
		\begin{algorithmic}[1]
			\\ \Lcomment{Input: initial vector $w_0\in \R^n$, number of iterations $J>0$}
			\\ \Lcomment{Input: Target accuracy $\epsilon>0$, failure probability $\delta \in(0,1)$}
			\State Sample to obtain approximate offsets: $\tilde{x}\in \R^{E}$ such that with probability $1-\frac{\delta}{2}$, $|\tilde{x}^{ab}_i-P^{ab}_i Lw_0|\leq \epsilon$ for all $i\in S$ and $a\in A_i$, $b\in B_{i,a}$:
			
			$\tilde{x}^{ab}_i=\texttt{ApxTransC}(w_0,\|L\|_\infty\|w_0\|_\infty,i,a,b,\epsilon,\frac{\delta}{2|E|})$
			\For{$j\in [J]$} $(w_j,\sigma_j,\tau_j)=\texttt{SApxVal} (w_{j-1},w_0,\tilde{x},\epsilon, \frac{\delta}{2J})$
			\EndFor
			\\ \Return $(w_J,\sigma_J,\tau_J)$
		\end{algorithmic}
	\end{small}
	\caption{Structured sampled randomized value iteration: \texttt{SSampledRandVI}{$(w_0,J,\epsilon,\delta)$}} \label{Algo5}
\end{algorithm}

\begin{algorithm}[h]
	\begin{small}
		\begin{algorithmic}[1]
			\\ \Lcomment{Input: Target accuracy $\epsilon>0$, failure probability $\delta \in(0,1)$}
			\State Let $K=\lceil\log_2(\frac{\|\psi\|_\infty W}{\epsilon})\rceil$ and $J= \lceil\frac{1}{1-\lambda}\log(4)\rceil$
			\State $w_0=0$ and $\epsilon_0=W$
			\For{$k\in[K]$} 
			\State $\epsilon_k=\frac {\epsilon_{k-1}}{2}=\frac{W}{2^k}$
			\State $\!\!(w_k,\sigma_k,\tau_k)\!=\!\texttt {SSampledRandVI}(w_{k-1},J,\frac{(1-\lambda)\epsilon_k}{4\|\psi^{-1}\|_\infty \Gamma}
                        ,\frac{\delta}{K})$
			\EndFor
			\\ \Return $(w_K,\sigma_K,\tau_K)$
		\end{algorithmic}
	\end{small}
	\caption{Structured sublinear randomized value iteration: \texttt{SSublinearRandVI}{$(\epsilon,\delta,\lambda,W,\Gamma,\|\psi^{-1}\|_\infty,\|\psi\|_\infty)$}}
	\label{Algo6}
\end{algorithm}

\begin{lemma}[adaptation of Lem.~4.2 in {\cite{sidford2018variance}}]\label{ApxTransLem} 
  \Cref{Algo1} runs in time $O(M^2\epsilon^{-2}\log(\frac{1}{\delta}))$, and it outputs $Y$ such that $|Y-P^{ab}_i u|\leq \epsilon$ with probability $1-\delta$.
\end{lemma}

\begin{lemma}[adaptation of Lem.~4.3 in {\cite{sidford2018variance}}]\label{ApxValComplexity}
	With probability $1-\delta$, \Cref{Algo2} returns $\tilde{w}$ such that $\|\tilde{w}-T(w)\|_\infty \leq 2 \Gamma \epsilon$, and then $\|\tilde{w}-T(w)\|_\psi \leq 2 \|\psi^{-1}\|_\infty \Gamma \epsilon$,
	and it runs in time:
	\[
	O\left(|E|\big\lceil\|w-w_0\|_\infty^2 \|L\|_\infty^2 \epsilon^{-2} \log\big(\frac{|E|}{\delta}\big)\big\rceil\right)\enspace.
	\]
\end{lemma}

\begin{lemma}\label{PsiContraction}
	If $w,w' \in \R^n$ satisfy $\|w'-T(w)\|_\psi\leq \alpha$ then $\|w'-w^\star\|_\psi\leq \alpha+\lambda \|w-w^\star\|_\psi$.
\end{lemma}

\begin{lemma}[adaptation of Lem.~4.5 in {\cite{sidford2018variance}}]\label{IneqRandVI}
	The sequence $(w_j)_{j\in [J]}$ generated by \Cref{Algo3} satisfies with probability $1-\delta$, that for all $j\in [J]$:
	\begin{equation*}
	\|w_j-w^\star\|_\psi \leq \frac{2\|\psi^{-1}\|_\infty \Gamma\epsilon}{1-\lambda}+\exp(-j(1-\lambda))\|w_0-w^\star\|_\psi
	\end{equation*}
	and if $J\geq \lceil\frac{1}{1-\lambda}\log(\frac{\|w_0-w^\star\|_\psi(1-\lambda)}{2 \|\psi^{-1}\|_\infty \Gamma \epsilon})\rceil$ then $\|w_J-w^\star\|_\psi \leq \frac{4 \|\psi^{-1}\|_\infty \Gamma\epsilon}{1-\lambda}$.
\end{lemma}

\begin{lemma}[adaptation of Lem.~4.6 in {\cite{sidford2018variance}}]\label{SRandVI}
	\Cref{Algo3} runs in time
	\begin{align*}
	O\bigg(|S\|E|+J|E|\bigg[\frac{\|\psi\|_\infty^2\|w_0-w^\star\|_\psi^2}{\epsilon^2}+\frac{\Gamma^2\|\psi\|_\infty^2\|\psi^{-1}\|_\infty^2}{(1-\lambda)^2} \bigg]\\
	\|L\|_\infty^2 \log\big(\frac{|E|J}{\delta}\big)\bigg)\enspace.
	\end{align*}
\end{lemma}

\begin{lemma}[adaptation of Lem.~4.8 and Lem.~4.9 in {\cite{sidford2018variance}}]\label{LinearComplexity}
	\Cref{Algo4} gives with probability $1-\delta$ that $\|w_k-w^\star\|_\psi\leq \epsilon_k$ for all $k\in[0,K]$, in particular $\|w_K-w^\star\|_\psi\leq \frac{\epsilon}{\|\psi\|_\infty}$ and then $\|w_K-w^\star\|_\infty\leq \epsilon$, and runs in time\footnote{As in \cite{sidford2018variance} we use $\tilde{O}$ to hide polylogarithmic factors in the input parameters, i.e. $\tilde{O} (f(x)) = O(f(x) \log(f(x))^{O(1)})$.}:
	\begin{align*}
	\tilde{O}\!\bigg(\!\!\big(|S\|E|\!+\!\frac{|E|\Gamma^2}{(1-\lambda)^3}\|\psi\|_\infty^2\|\psi^{-1}\|_\infty^2 \big)
	\|L\|_\infty^2\!\log(\frac{W}{\epsilon})\log(\frac{1}{\delta})\!\!\bigg).
	\end{align*}
\end{lemma}

\begin{lemma}[adaptation of Lem.~4.10 and Lem.~4.12 in {\cite{sidford2018variance}}]\label{SublinearComplexity}
	\Cref{Algo6} gives with probability $1-\delta$ that $\|w_k-w^\star\|_\psi\leq \epsilon_k$ for all $k\in[0,K]$, in particular $\|w_K-w^\star\|_\psi\leq \frac{\epsilon}{\|\psi\|_\infty}$ and then $\|w_K-w^\star\|_\infty\leq \epsilon$, and runs in time
	\begin{align*}
	\tilde{O}\bigg(|E|\Gamma^2 \|\psi\|_\infty^2\|\psi^{-1}\|_\infty^2 \bigg[\frac{\|\psi\|_\infty^2 W^2}{(1-\lambda)^2\epsilon^2}+\frac{1}{(1-\lambda)^3}\bigg]\\
	\|L\|_\infty^2\log(\frac{1}{\delta})\bigg)\enspace.
	\end{align*}
\end{lemma}

\section{Variance reduced deflated value iteration for ergodic problems}\label{sec-ergodic}
To solve the mean-payoff problem of~\cref{MPP}, we consider the equation:
\begin{equation}\label{MPSG}
\eta e+v=T(v) \; \text{and} \; v_c=0,\; \eta\in \R,\;v\in \R^n
\end{equation}
where $T$ is as in \Cref{def-shapley}: $T_i(v)=\min_{a\in A_i} \max_{b\in B_{i,a}}\big\{r_i^{ab}+ \sum_{j\in S}{P_{ij}^{ab}} v_j\big\}, \forall i\in S$.
Throughout the section, we make Assumption~\ref{assum-commonc}.
We denote by $(\eta^\star,v^\star)$ the unique solution of this problem. 
We know by \Cref{RhoHypothesis} that there exists a vector
$\myphi \in \R_{+}^n$ satisfying \eqref{PhiIneq}.

\Cref{ErgToDisc} shows that \eqref{MPSG} is equivalent to the equation:
\begin{equation}\label{TPhi}
T^{\myphi}(w)=w,\quad w\in \R^n \enspace,
\end{equation}
with $\eta=w_c$, $v=\myphi(w-w_c)$, and $T^{\myphi}$ is given by \eqref{T-Phi}.

\subsection{Computing an h-transform of the ergodic problem}
Here we want to find a vector $\myphi \in \R_{+}^n$ satisfying \eqref{PhiIneq}.
First, we consider the problem of finding the vector of maximal expected first hitting times of state $c$, denoted $\myphi^\star$ as in \eqref{PhiTmax},
and we suppose that we know a bound $H$ on it: 
\begin{equation}\label{Hittingtime}
H\geq \|\myphi^\star\|_\infty=\max_{i\in S}\mc{T}_{ic}\enspace.
\end{equation}
We define the scalar $\lambda\in [0,1)$ by
\begin{equation}\label{lambda}
\lambda:=1-1/H\geq 1-1/\|\myphi^\star\|_\infty\enspace.
\end{equation}

We remark that the component $\myphi^\star_c$ can be computed in time $O(|E|)$ from the other components since $\myphi^\star_c=1+\max_{a\in A_i,b\in B_{i,a}}[\sum_{j \in S, j\neq c}P^{ab}_{ij}\myphi^\star_j]$. By considering $w^\star=(\myphi^\star_i)_{i\in S\setminus\{c\}}\in \R^{n-1}$ and the matrices $\tilde{P}^{\sigma\tau}\in \R^{(n-1)\times (n-1)}$ defined from $P^{\sigma\tau}$ by removing the $c$ row and the $c$ column, the problem becomes 
\begin{equation}\label{Tm}
w^\star=T^m(w^\star) \enspace,
\end{equation}
where the operator $T^m$ is such that
\begin{equation*}
T^m_i(w)=1+\max_{a,b}[\tilde{P}^{ab}_{i}w], \quad \forall i\in S\setminus\{c\}, \forall w\in \R^{n-1}
\end{equation*}

The operator $T^m$ is a particular case of the operator $T$ in (\ref{T}): 
for $(\sigma,\tau) \in \polmin\times\polmax$, $P^{\sigma \tau}=\tilde{P}^{\sigma \tau}$ are sub-Markovian. For $(i,a,b)\in E$, $\gamma^{ab}_{i}=1$ and then $\Gamma=1$, $G^{ab}_{i}(w) =1$ for all $w\in \R^{n-1}$ and $L$ is the identity then $\|L\|_\infty=1$.

From \Cref{PhiContraction} and \Cref{lambda}, we know that the operator $T^m$ is $\mu-$contracting in the sup-norm $\|\cdot\|_{w^\star}$, with $\mu:=1-1/\|w^\star\|_\infty\leq 1-1/\|\myphi^\star\|_\infty\leq \lambda$, so here we take $\psi=w^\star$, and $\lambda$ as contraction rate.
We have $\|w^{\star-1}\|_\infty\leq 1$ and $\|w^\star\|_\infty\leq \|\myphi^\star\|_\infty\leq \frac{1}{1-\lambda}$, then according to \Cref{d1d2} we can take in \Cref{Algo4} and \Cref{Algo6} the scalars $1$ and $\frac{1}{1-\lambda}$ instead of $\|w^{\star-1}\|_\infty$ and $\|w^\star\|_\infty$ respectively.
We have $\|w^\star\|_{w^\star}= 1$, then we take $W=1$.

We use \Cref{Algo4} and \Cref{Algo6} to find an $\epsilon-$approximation of $w^\star$ in near linear and sublinear time respectively. Then we can deduce $\myphi'$ an $\epsilon-$approximation of $\myphi^\star$. By taking $\epsilon=\frac{1}{4}$ and considering $\myphi=2\myphi'$, we deduce the following \Cref{PhiLinear} and \Cref{PhiSubLinear}. In both theorems $\lambda_\myphi:=1-1/\|\myphi\|_\infty$ satisfies $\frac{1}{1-\lambda_\myphi}\leq\frac{2}{1-\lambda}+\frac{1}{2}$.

\begin{theorem}\label{PhiLinear}
	By calling \Cref{Algo4}, we can find $\myphi \in \R_{+}^n$ such that $\myphi_i\geq 1+\max_{a,b}[P^{ab}_{(c)i}\myphi]$ for all $i\in S$, in time
	\begin{equation*}
	\tilde{O}\bigg(\bigg(|S\|E|+\frac{|E|}{(1-\lambda)^5}\bigg)\log(\frac{1}{\delta})\bigg)\enspace.
	\end{equation*}
\end{theorem}

\begin{theorem}\label{PhiSubLinear}
	By calling \Cref{Algo6}, we can find $\myphi \in \R_{+}^n$ such that $\myphi_i\geq 1+\max_{a,b}[P^{ab}_{(c)i}\myphi]$ for all $i\in S$ in time
	\begin{equation*}
	\tilde{O}\bigg(\frac{|E|}{(1-\lambda)^6}\log(\frac{1}{\delta})\bigg)\enspace.
	\end{equation*}
\end{theorem}

\subsection{Solving the ergodic problem}
We suppose that we have identified a vector $\myphi$ satisfying \eqref{PhiIneq}. Now we consider the equation \eqref{TPhi}.

As in \cref{sec-SSG}, we can write $T^{\myphi}$ in the general form given by \Cref{Tstructured},
where $\gamma^{ab}_{i}=\myphi_i^{-1}$ for all $i\in S, a\in A_i, b\in B_{i,a}$ and then $\Gamma=\max_{(i,a,b)\in E}\gamma^{ab}_{i}\leq1$, $L:w\mapsto \myphi(w-w_ce)$ is a sparse linear operator and we have $\|L\|_\infty\leq 2\|\myphi\|_\infty=\frac{2}{1-\lambda_\myphi}$ and $G^{ab}_{i}:w \mapsto \myphi_i^{-1}r^{ab}_{i}+w_c(1-\myphi_i^{-1})$ is also a sparse affine operator for all $i\in S, a\in A_i, b\in B_{i,a}$.

By \Cref{TPhiContraction}, $T^{\myphi}$ is a $\lambda_\myphi-$contraction in the sup-norm $\|\cdot\|_\infty$, where $\lambda_\myphi= 1-1/\|\myphi\|_\infty$. We take here $\psi=e$, which means $\|\cdot\|_\psi=\|\cdot\|_\infty$ and $\|\psi^{-1}\|_\infty=\|\psi\|_\infty=1$.

By \Cref{boundPay}, we have $\|w^\star\|_\infty\leq R$, where $w^\star$ is the solution of \Cref{TPhi}. So we take $W=R$.

The following theorems give an approximation for $w^\star$ and then also for $v^\star$ and $\eta^\star$ in nearly linear time with \Cref{Linear} (based on \Cref{PhiLinear} and \Cref{LinearComplexity}), and in sublinear time with \Cref{SubLinear} (based on \Cref{PhiSubLinear} and \Cref{SublinearComplexity}). The time complexities considered include the times needed to find $\myphi$.

\begin{theorem}\label{Linear}
	With probability $1-\delta$, we find $\myphi$ satisfying \eqref{PhiIneq} and the call of \Cref{Algo4},  \texttt{SHighPrecisionRandVI}$(\epsilon,\frac{\delta}{2},\lambda_\myphi,R,1,1,1)$ returns $w\in \R^n$ such that $\|w-w^\star\|_\infty\leq \epsilon$. Therefore we obtain $\eta=w_{c}$ and $v=\myphi(w-w_{c}e)$ such that $|\eta-\eta^\star|\leq \epsilon$ and $\|v-v^\star\|_\infty\leq \frac{5\epsilon}{1-\lambda}$. The run time needed is
	\begin{equation*}
	\tilde{O}\bigg(\bigg(|S\|E|+\frac{|E|}{(1-\lambda)^5}\bigg)\log(\frac{R}{\epsilon})\log(\frac{1}{\delta})\bigg)\enspace.
	\end{equation*}
\end{theorem}

\begin{theorem}\label{SubLinear}
	With probability $1-\delta$, we find $\myphi$ satisfying \eqref{PhiIneq} and the call of \Cref{Algo6},  \texttt{SSublinearRandVI}$(\epsilon,\frac{\delta}{2},\lambda_\myphi,R,1,1,1)$ returns $w\in \R^n$ such that $\|w-w^\star\|_\infty\leq \epsilon$. Therefore we obtain $\eta=w_{c}$ and $v=\myphi(w-w_{c}e)$ such that $|\eta-\eta^\star|\leq \epsilon$ and $\|v-v^\star\|_\infty\leq \frac{5\epsilon}{1-\lambda}$. The run time needed is
	\begin{equation*}
	\tilde{O}\bigg(|E|\bigg[\frac{R^2}{(1-\lambda)^4\epsilon^2}+\frac{1}{(1-\lambda)^6}\bigg]\log(\frac{1}{\delta})\bigg)\enspace.
	\end{equation*}
\end{theorem}

\section{Comparison with alternative approaches}
\label{sec-example}
\subsection{A cyclic example: no contraction or mixing}
The convergence proof of the sampled relative value iteration
in \cite{jain} requires the Dobrushin ergodicity coefficient
\begin{align}
  \alpha = 1- \min_{i,j\in[n],a\in A_i, a'\in A_j } \sum_{k\in[n]}\min(P^a_{ik},P^{a'}_{jk})
  \label{e-dobrushin}
\end{align}
to be smaller than $1$.
Consider the $0$-player instance
with a cyclic matrix $P$, presented in \Cref{ex-illustrate}.
Here, $\alpha=1$, and actually,
relative value iteration does not converge. 
Moreover, the mixing time used in the bound of \cite{wang2017primal} is infinite.
However, as shown in~\Cref{ex-illustrate}, the deflation+h-transform
methods reduces to a fixed point problem with a contraction
rate of $1/2$.

\subsection{An example with small hitting times}
\label{sec-smallhitting}
Let us first consider a $0-$player problem with state space $[n]$, $T(v)=r+Qv, \forall v\in \R^n$, where $r\in \R^n$ fixed and the probability transition matrix $Q\in \R^{n\times n}$, 
such that
\begin{align*}
Q_{i,1}=Q_{i,i+1}=1/2, i\in[n-1],\quad Q_{n,1}=1\enspace .
\end{align*}
We can easily prove that the expected first return time to state $n$  is $\mc{T}_{cc}=\Omega(2^n)$. If we denote by $\nu$ the stationary distribution of $Q$,
it follows that $\nu_n=O(2^{-n})$. In \cite{wang2017primal}, one supposes
the stationary distribution of $Q$ satisfies $\frac{1}{\sqrt{\tau}n}e\leq \nu\leq \frac{\sqrt{\tau}}{n}e$, so, $\tau= \Omega(2^{2n}/n^2)$.
The complexity bound of \cite{wang2017primal} is exponential in this example,
since it includes a $\tau^2$ factor. 

By using our technique, we will first choose $c=1$ and we verify easily that the first hitting time vector $\myphi^\star$ satisfies $\|\myphi^\star\|_\infty\leq 2$ (more precisely $\myphi^\star_i=2-\frac{1}{2^{n-i}}, \forall i\in[n]$). Then the new operator $T^{\myphi^\star}$ is $1/2-$contracting which leads to fast convergence. In particular, \Cref{Linear} gives a time complexity $\tilde{O}(n^2\log(\frac{R}{\epsilon})\log(\frac{1}{\delta}))$, and \Cref{SubLinear} gives a time complexity $\tilde{O}(n\frac{R^2}{\epsilon^2}\log(\frac{1}{\delta}))$, where $R$ is an upper bound on the payments, $\epsilon$ is the target accuracy and $\delta$ is the failure probability.

In this $0$-player example, $\alpha=1/2$, so we could use relative value iteration. This is no longer the case if we consider the following
$1$-player variant. By identifying $n+1$ and $1$, consider the stochastic matrix
\begin{align*}
Q'_{i+1,2}=Q'_{i+1,i+2}=\frac 1 2, i\in[n-1],\quad Q'_{1,2}=1\enspace .
\end{align*}
let $r'\in\R^n$ be a vector of payments, and consider the Bellman operator
\(
T(x) =\max(r+Qx, r'+Q'x) \),
so that there are two actions in every state.
Now, $\alpha=1$ and the convergence of the relative value iteration~\cite{jain} is not guaranteed. However, we observe that for all policies, the probability
to reach the set of states $\{1,2\}$ in one step is at least $1/2$.  Moreover,
for all actions, the probability of the transition $1\to 2$
is
at least
$1/2$.
It follows that by choosing $c=2$,
the vector of maximal hitting times $\myphi^\star$ satisfies
$\|\myphi^\star\|_\infty \leq 4$
(more precisely $\myphi^\star_1=2$ and $\myphi^\star_i=4-\frac{1}{2^{n-i}}, \forall i\in\{2,\cdots,n\}$). Hence, the deflation+h-transform method
still has a sublinear behavior on this example.

\nocite{1806.01492}

{\em Acknowledgement}.\/ We thank the referees for detailed and helpful remarks.

\bibliographystyle{IEEEtran}

\bibliography{variancereductionergodic}

\begin{appendices}

\end{appendices}
\end{document}